\documentclass[11pt,leqno]{article}

\usepackage{theorem,amssymb,amstex}

\makeatletter
\ifx\@auxdefsloaded\relax\endinput
\else\let\@auxdefsloaded\relax\fi

\def\newenvironment{%
   \@ifnextchar *{\@@newenv{\global\@ignoretrue}}{\@@newenv{}*}}

\def\@@newenv#1*#2{%
   \@ifnextchar [{\@newenv{#1}{#2}}{\@newenv{#1}{#2}[0]}}

\long\def\@newenv#1#2[#3]#4#5{%
   \expandafter\newcommand\csname#2\endcsname[#3]{#4}%
   \expandafter\long\expandafter\def\csname end#2\endcsname{#5#1}}

\def\renewenvironment{%
   \@ifnextchar *{\@@renewenv{\global\@ignoretrue}}{\@@renewenv{}*}}

\def\@@renewenv#1*#2{%
   \@ifnextchar [{\@renewenv{#1}{#2}}{\@renewenv{#1}{#2}[0]}}

\long\def\@renewenv#1#2[#3]#4#5{%
   \expandafter\renewcommand\csname#2\endcsname[#3]{#4}%
   \expandafter\long\expandafter\def\csname end#2\endcsname{#5#1}}

\def\newoptcommand#1#2{%
   \@ifnextchar [{\@optargdef#1#2}{\@optargdef#1#2[1]}}

\def\renewoptcommand#1#2{%
   \edef\@tempa{\expandafter\@cdr\string#1\@nil}%
   \@ifundefined{\@tempa}{%
      \@latexerr{\string#1\space undefined}\@ehc}{}%
   \@ifnextchar [{\@reoptargdef#1#2}{\@reoptargdef#1#2[1]}}

\long\def\@optargdef#1#2[#3]#4{%
   \@ifdefinable #1{\@reoptargdef#1#2[#3]{#4}}}

\long\def\@reoptargdef#1#2[#3]#4{%
   \@tempcnta#3\relax \@tempcntb \@ne
   \let#1\relax \let\@tempa\relax
   \edef\@tempb{\long\def\csname\string#1\endcsname
      [\@tempa\the\@tempcntb]}%
   \advance\@tempcntb \@ne \advance\@tempcnta \m@ne
   \@whilenum\@tempcnta>0\do{%
      \edef\@tempb{\@tempb\@tempa\the\@tempcntb}%
      \advance\@tempcntb \@ne \advance\@tempcnta \m@ne}%
   \let\@tempa=##\@tempb{#4}%
   \def#1{\@ifnextchar [{\csname\string#1\endcsname}{%
      \csname\string#1\endcsname[#2]}}}

\def\newoptenvironment{%
   \@ifnextchar *{\@@newoptenv{\global\@ignoretrue}}{%
      \@@newoptenv{}*}}

\def\@@newoptenv#1*#2#3{%
   \@ifnextchar [{\@newoptenv{#1}{#2}{#3}}{%
      \@newoptenv{#1}{#2}{#3}[0]}}

\long\def\@newoptenv#1#2#3[#4]#5#6{%
   \expandafter\newoptcommand\csname#2\endcsname{#3}[#4]{#5}%
   \expandafter\long\expandafter\def\csname end#2\endcsname{#6#1}}

\def\renewoptenvironment{%
   \@ifnextchar *{\@@renewoptenv{\global\@ignoretrue}}{%
      \@@renewoptenv{}*}}

\def\@@renewoptenv#1*#2#3{%
   \@ifnextchar [{\@renewoptenv{#1}{#2}{#3}}{%
      \@renewoptenv{#1}{#2}{#3}[0]}}

\long\def\@renewoptenv#1#2#3[#4]#5#6{%
   \expandafter\renewoptcommand\csname#2\endcsname{#3}[#4]{#5}%
   \expandafter\long\expandafter\def\csname end#2\endcsname{#6#1}}

\newcounter{keepoptional}
\newcounter{optnestctr}

\newoptenvironment*{optional}{1}[1]{%
   \ifnum#1>\value{keepoptional}\relax
      \setcounter{optnestctr}{0}\@powerup\expandafter\@endopt\fi
   \ignorespaces}{}

\def\@powerup{\catcode`\{=12 \catcode`\}=12 \catcode`\\=12 \relax}
\def\@powerdown{\catcode`\{=1 \catcode`\}=2 \catcode`\\=0 \relax}
\begingroup \catcode`|=0 \catcode`[=1 \catcode`]=2 \@powerup
   |long|gdef|@endopt#1\end{optional}[%
      |@beginopt#1\begin{optional}*]%
   |long|gdef|@beginopt#1\begin{optional}[%
      |@ifstar[%
         |ifnum|value[optnestctr]>0|relax
            |addtocounter[optnestctr][-1]|let|@zzzap=|@endopt
         |else
            |@powerdown|end[optional]|let|@zzzap=|relax|fi
         |@zzzap]%
      [
         |addtocounter[optnestctr][1]|@beginopt]]%
|endgroup

\makeatother
\makeatletter
\ifx\@auxdefsloaded\relax
\else \input{auxdefs.sty}\fi
\newskip\dgARROWLENGTH  \dgARROWLENGTH=2.5em\relax
\newskip\dgHORIZPAD     \dgHORIZPAD=1em\relax
\newskip\dgVERTPAD      \dgVERTPAD=2ex\relax
\newskip\dgLABELOFFSET  \dgLABELOFFSET=.7ex\relax
\newcount\dgARROWPARTS  \dgARROWPARTS=4\relax
\newcommand{\dgeverynode}{\displaystyle}
\newcommand{\dgeverylabel}{\scriptstyle}
\newskip\dgDOTSPACING   \dgDOTSPACING=0.35em
\newskip\dgDOTSIZE      \dgDOTSIZE=1.5\fontdimen8\tenln
\newcount\dgMAXSQUARE   \dgMAXSQUARE=4\relax
\newcount\dgMINDBLSQ    \dgMINDBLSQ=10\relax
\newskip\dgCOLUMNWIDTH  \dgCOLUMNWIDTH=2em\relax
\chardef\f@ur=4
\@namedef{dgo@}{\let\dg@VECTOR=\vector
   \dg@LBLPOS=\dgARROWPARTS \divide\dg@LBLPOS\tw@}%
\@namedef{dgo@-}{\let\dg@VECTOR=\line}%
\@namedef{dgo@!}{\let\dg@VECTOR=\dg@novector}%
\@namedef{dgo@1}{\dg@LBLPOS=\@ne\relax}%
\@namedef{dgo@2}{\dg@LBLPOS=\tw@\relax}%
\@namedef{dgo@3}{\dg@LBLPOS=\thr@@\relax}%
\@namedef{dgo@4}{\dg@LBLPOS=\f@ur\relax}%
\@namedef{dgo@5}{\dg@LBLPOS=5\relax}%
\@namedef{dgo@6}{\dg@LBLPOS=6\relax}%
\@namedef{dgo@7}{\dg@LBLPOS=7\relax}%
\@namedef{dgo@8}{\dg@LBLPOS=8\relax}%
\@namedef{dgo@9}{\dg@LBLPOS=9\relax}%
\def\dgt@e{\dg@DX=\@ne \dg@DY=\z@ \dg@SIZE=\@ne}%
\def\dgt@w{\dg@DX=\m@ne \dg@DY=\z@ \dg@SIZE=\@ne}%
\def\dgt@n{\dg@DX=\z@ \dg@DY=\@ne \dg@SIZE=\@ne}%
\def\dgt@s{\dg@DX=\z@ \dg@DY=\m@ne \dg@SIZE=\@ne}%
\def\dgt@ne{\dg@DX=\@ne \dg@DY=\@ne \dg@SIZE=\@ne}%
\def\dgt@se{\dg@DX=\@ne \dg@DY=\m@ne \dg@SIZE=\@ne}%
\def\dgt@nw{\dg@DX=\m@ne \dg@DY=\@ne \dg@SIZE=\@ne}%
\def\dgt@sw{\dg@DX=\m@ne \dg@DY=\m@ne \dg@SIZE=\@ne}%
\def\dgt@nne{\dg@DX=\@ne \dg@DY=\tw@ \dg@SIZE=\@ne}%
\def\dgt@nnw{\dg@DX=\m@ne \dg@DY=\tw@ \dg@SIZE=\@ne}%
\def\dgt@sse{\dg@DX=\@ne \dg@DY=-\tw@ \dg@SIZE=\@ne}%
\def\dgt@ssw{\dg@DX=\m@ne \dg@DY=-\tw@ \dg@SIZE=\@ne}%
\def\dgt@ene{\dg@DX=\tw@ \dg@DY=\@ne \dg@SIZE=\tw@}%
\def\dgt@ese{\dg@DX=\tw@ \dg@DY=\m@ne \dg@SIZE=\tw@}%
\def\dgt@wnw{\dg@DX=-\tw@ \dg@DY=\@ne \dg@SIZE=\tw@}%
\def\dgt@wsw{\dg@DX=-\tw@ \dg@DY=\m@ne \dg@SIZE=\tw@}%
\def\dggeometry{
   \dg@ZTEMP=\dg@XGRID \multiply\dg@ZTEMP\tw@
   \ifnum\dg@YGRID=\z@ \dg@ZTEMP=\tw@
   \else \divide\dg@ZTEMP\dg@YGRID \fi
   \ifnum\dg@ZTEMP>\f@ur \dg@ZTEMP=\f@ur \fi
   \ifnum\dg@ZTEMP<\@ne \dg@ZTEMP=\@ne \fi
   \unitlength=2sp\relax
   \ifnum\dg@ZTEMP<\tw@
      \advance\dg@ZTEMP\@ne
      \multiply\unitlength\dg@YGRID
   \else
      \multiply\unitlength\dg@XGRID \divide\unitlength\dg@ZTEMP
   \fi
   \dg@XGRID=\dg@ZTEMP \dg@YGRID=\tw@
   \dg@rmcommondiv\tw@\dg@XGRID\dg@YGRID
   \divide\unitlength\dg@YGRID \divide\unitlength\@m\relax}
\@namedef{dgo@..}{\let\dg@VECTOR=\dg@dotvector}%

\def\dg@dotvector(#1,#2)#3{%
   \begingroup
   \dg@XTEMP=#1\relax \dg@YTEMP=#2\relax
   \let\dg@NDOTS=\dg@XEND \let\dg@DOTDIAM=\dg@WEND
   \dg@NDOTS=\unitlength \multiply\dg@NDOTS #3\relax
   \dg@ZTEMP=\dg@YTEMP \dg@changesign\dg@YTEMP\dg@ZTEMP
   \ifnum\dg@XTEMP>\z@
      \ifnum\dg@YTEMP>\dg@XTEMP
         \multiply\dg@NDOTS\dg@YTEMP \divide\dg@NDOTS\dg@XTEMP \fi
   \else\ifnum\dg@XTEMP<\z@
      \ifnum\dg@YTEMP>-\dg@XTEMP
         \multiply\dg@NDOTS\dg@YTEMP \divide\dg@NDOTS-\dg@XTEMP \fi
   \fi\fi
   \dg@YTEMP=\dg@ZTEMP
   \divide\dg@NDOTS\dgDOTSPACING
   \ifnum\dg@NDOTS>\z@\else \dg@NDOTS=\@ne \fi
   \dg@ZTEMP=\unitlength \multiply\dg@ZTEMP #3\relax
   \divide\dg@ZTEMP\dg@NDOTS
   \ifnum\dg@XTEMP=\z@
      \dg@changesign\dg@ZTEMP\dg@YTEMP \dg@YTEMP=\dg@ZTEMP
   \else
      \dg@changesign\dg@ZTEMP\dg@XTEMP
      \multiply\dg@YTEMP\dg@ZTEMP \divide\dg@YTEMP\dg@XTEMP
      \dg@XTEMP=\dg@ZTEMP
   \fi
   \divide\dg@XTEMP\unitlength \divide\dg@YTEMP\unitlength
   \begin{picture}(0,0)
      \dg@DOTDIAM=\dgDOTSIZE \divide\dg@DOTDIAM\unitlength
      \multiput(0,0)(\dg@XTEMP,\dg@YTEMP){\dg@NDOTS}{%
         \circle*{\dg@DOTDIAM}}%
      \multiply\dg@XTEMP\dg@NDOTS \multiply\dg@YTEMP\dg@NDOTS
      \put(\dg@XTEMP,\dg@YTEMP){\vector(#1,#2){0}}%
   \end{picture}%
   \endgroup}%
\newif\ifdg@LATEXGEOM \dg@LATEXGEOMfalse

\@ifundefined{lamsvector}{}{%
   \@namedef{dgo@}{%
      \let\dg@VECTOR=\lamsvector
      \lamsshaft{?}\lamssource{?}\lamstarget{?}%
      \dg@LBLPOS=\dgARROWPARTS \divide\dg@LBLPOS\tw@\relax}%
   \@namedef{dgo@..}{\lamsshaft{.}}%
   \@namedef{dgo@=}{\lamsshaft{=}}%
   \@namedef{dgo@-}{\lamstarget{0}}%
   \@namedef{dgo@'}{\lamstarget{'}}%
   \@namedef{dgo@`}{\lamstarget{`}}%
   \@namedef{dgo@A}{\lamstarget{A}}%
   \@namedef{dgo@V}{\lamssource{V}}%
   \@namedef{dgo@J}{\lamssource{J}}%
   \@namedef{dgo@L}{\lamssource{L}}%
   \@namedef{dgo@S}{\lamssource{S}}%
   \def\dggeometry{
      \dg@ZTEMP=\dg@XGRID \multiply\dg@ZTEMP\tw@
      \ifnum\dg@YGRID=\z@ \dg@ZTEMP=\tw@
      \else \divide\dg@ZTEMP\dg@YGRID \fi
      \ifnum\dg@ZTEMP>6\relax \dg@ZTEMP=6\relax \fi
      \ifdg@LATEXGEOM\ifnum\dg@ZTEMP>\f@ur \dg@ZTEMP=\f@ur \fi\fi
      \ifnum\dg@ZTEMP<\@ne \dg@ZTEMP=\@ne \fi
      \unitlength=2sp\relax
      \ifnum\dg@ZTEMP<\tw@
         \advance\dg@ZTEMP\@ne
         \multiply\unitlength\dg@YGRID
      \else
         \multiply\unitlength\dg@XGRID \divide\unitlength\dg@ZTEMP
      \fi
      \dg@XGRID=\dg@ZTEMP \dg@YGRID=\tw@
      \dg@rmcommondiv\tw@\dg@XGRID\dg@YGRID
      \divide\unitlength\dg@YGRID \divide\unitlength\@m
      \dg@LATEXGEOMfalse}
   \def\dgt@nee{\dg@DX=\tw@ \dg@DY=\@ne \dg@SIZE=\tw@}%
   \def\dgt@see{\dg@DX=\tw@ \dg@DY=\m@ne \dg@SIZE=\tw@}%
   \def\dgt@nww{\dg@DX=-\tw@ \dg@DY=\@ne \dg@SIZE=\tw@}%
   \def\dgt@sww{\dg@DX=-\tw@ \dg@DY=\m@ne \dg@SIZE=\tw@}%
   \def\dgt@nnne{\dg@DX=\@ne \dg@DY=\thr@@ \dg@SIZE=\@ne}%
   \def\dgt@nnnw{\dg@DX=\m@ne \dg@DY=\thr@@ \dg@SIZE=\@ne}%
   \def\dgt@sssw{\dg@DX=\m@ne \dg@DY=-\thr@@ \dg@SIZE=\@ne}%
   \def\dgt@ssse{\dg@DX=\@ne \dg@DY=-\thr@@ \dg@SIZE=\@ne}%
   \def\dgt@nnnee{\dg@DX=\tw@ \dg@DY=\thr@@ \dg@SIZE=\tw@}%
   \def\dgt@nnnww{\dg@DX=-\tw@ \dg@DY=\thr@@ \dg@SIZE=\tw@}%
   \def\dgt@sssww{\dg@DX=-\tw@ \dg@DY=-\thr@@ \dg@SIZE=\tw@}%
   \def\dgt@sssee{\dg@DX=\tw@ \dg@DY=-\thr@@ \dg@SIZE=\tw@}%
   \def\dgt@nneee{\dg@DX=\thr@@ \dg@DY=\tw@ \dg@SIZE=\thr@@}%
   \def\dgt@nnwww{\dg@DX=-\thr@@ \dg@DY=\tw@ \dg@SIZE=\thr@@}%
   \def\dgt@sswww{\dg@DX=-\thr@@ \dg@DY=-\tw@ \dg@SIZE=\thr@@}%
   \def\dgt@sseee{\dg@DX=\thr@@ \dg@DY=-\tw@ \dg@SIZE=\thr@@}%
   \def\dgt@neee{\dg@DX=\thr@@ \dg@DY=\@ne \dg@SIZE=\thr@@
      \global\dg@LATEXGEOMtrue}%
   \def\dgt@nwww{\dg@DX=-\thr@@ \dg@DY=\@ne \dg@SIZE=\thr@@
      \global\dg@LATEXGEOMtrue}%
   \def\dgt@swww{\dg@DX=-\thr@@ \dg@DY=\m@ne \dg@SIZE=\thr@@
      \global\dg@LATEXGEOMtrue}%
   \def\dgt@seee{\dg@DX=\thr@@ \dg@DY=\m@ne \dg@SIZE=\thr@@
      \global\dg@LATEXGEOMtrue}%
}
\newcount\dg@HORIZ      \newcount\dg@VERT
\newcount\dg@XLPAD      \newcount\dg@YBPAD
\newcount\dg@XRPAD      \newcount\dg@YTPAD
\newcount\dg@X          \newcount\dg@Y
\newcount\dg@XGRID      \newcount\dg@YGRID
\newcount\dg@SIZE
\newcount\dg@USERSIZE
\newcount\dg@DX         \newcount\dg@DY
\newcount\dg@XLBL       \newcount\dg@YLBL
\newcount\dg@XOFFSET    \newcount\dg@YOFFSET
\newcount\dg@LBLOFF
\newcount\dg@XLBLOFF    \newcount\dg@YLBLOFF
\newcount\dg@LBLPOS
\newcount\dg@XTEMP      \newcount\dg@YTEMP
\newcount\dg@ZTEMP
\newcount\dg@XNODE      \newcount\dg@YNODE
\newcount\dg@XEND       \newcount\dg@YEND
\newcount\dg@WEND       \newcount\dg@HEND
\newcount\dg@COUNT
\newbox\dg@NODEBOX
\newoptenvironment*{diagram}{}[1]{%
   \global\advance\dg@COUNT\@ne \typeout{[diagram \the\dg@COUNT:}%
   \let\node=\dg@node \let\\=\dg@cr \let\arrow=\dg@arrow
   \def\dg@BIGNODE{#1}%
   \ifx\dg@BIGNODE\@empty\else
      \setbox\dg@NODEBOX=\hbox{$\dgeverynode{#1}$}%
      \dg@XTEMP=\wd\dg@NODEBOX
      \dg@YTEMP=\ht\dg@NODEBOX \advance\dg@YTEMP\dp\dg@NODEBOX
      \ifnum\dg@YTEMP=\z@ \dg@ZTEMP=\z@
      \else \dg@ZTEMP=\dg@XTEMP \divide\dg@ZTEMP\dg@YTEMP\fi
      \ifnum\dg@ZTEMP>\dgMAXSQUARE
         \ifnum\dg@ZTEMP<\dgMINDBLSQ \dg@XGRID=\thr@@ \dg@YGRID=\tw@
         \else \dg@XGRID=\tw@ \dg@YGRID=\@ne \fi
      \else \dg@XGRID=\@ne \dg@YGRID=\@ne \fi
      \advance\dg@XTEMP\dgHORIZPAD \advance\dg@YTEMP\dgVERTPAD
      \dg@XLPAD=\dg@XTEMP \divide\dg@XLPAD\tw@
      \dg@XRPAD=\dg@XTEMP \divide\dg@XRPAD\tw@
      \dg@YBPAD=\dg@YTEMP \divide\dg@YBPAD\tw@
      \dg@YTPAD=\dg@YTEMP \divide\dg@YTPAD\tw@
      \advance\dg@XTEMP\dgARROWLENGTH
      \divide\dg@XTEMP\dg@XGRID \divide\dg@XTEMP\@m
      \unitlength=1sp\relax \multiply\unitlength\dg@XTEMP
   \fi
   \typeout{saving...}%
   \dg@X=-\@ne \dg@Y=\z@ \dg@HORIZ=\z@\relax
   \let\dg@SLIST=\@empty
   \let\dg@NLIST=\@empty \let\dg@ALIST=\@empty
   \let\dg@PASS=\dg@savepass
}{
   \dg@VERT=-\dg@Y\relax
   \typeout{calculating...}%
   \ifx\dg@BIGNODE\@empty
      \dg@XGRID=\z@ \dg@YGRID=\z@
      \dg@XLPAD=\z@ \dg@XRPAD=\z@ \dg@YBPAD=\z@ \dg@YTPAD=\z@
      \let\dg@PASS=\dg@geompass
      \dg@NLIST \dg@ALIST
      \ifnum\dg@XGRID>\z@\else \dg@XGRID=\quad\relax \fi
      \ifnum\dg@YGRID>\z@\else \dg@YGRID=\quad\relax \fi
      \dggeometry
      \ifdim\unitlength>\z@\else \unitlength=\quad \fi
      \ifnum\dg@XGRID>\z@\else \dg@XGRID=\@ne \fi
      \ifnum\dg@YGRID>\z@\else \dg@YGRID=\@ne \fi
   \fi
   \dg@LBLOFF=\dgLABELOFFSET \divide\dg@LBLOFF\unitlength
   \multiply\dg@HORIZ\@m \multiply\dg@HORIZ\dg@XGRID
   \multiply\dg@VERT\@m \multiply\dg@VERT\dg@YGRID
   \divide\dg@XLPAD\unitlength \divide\dg@XRPAD\unitlength
   \divide\dg@YBPAD\unitlength \divide\dg@YTPAD\unitlength
   \typeout{drawing...}%
   \let\dg@PASS=\dg@drawpass
   \hspace*{\dg@XLPAD\unitlength}%
   \vcenter{%
      \hsize=0pt\relax\parindent=0pt\relax
      \parskip=0pt\relax\baselineskip=0pt\relax
      \vspace*{\dg@YTPAD\unitlength}%
      \begin{picture}(0,0)(0,0)\dg@NLIST\dg@ALIST\end{picture}%
      \raisebox{\z@}[\z@][\dg@VERT\unitlength]{}%
      \vspace*{\dg@YBPAD\unitlength}%
   }
   \hspace*{\dg@HORIZ\unitlength}%
   \hspace*{\dg@XRPAD\unitlength}%
   \typeout{done]}}%

\def\dg@savepass{s}
\def\dg@geompass{g}
\def\dg@drawpass{d}

\newoptcommand{\dg@node}{\@ne}[2]{%
   \ifx\dg@PASS\dg@savepass
      \dg@XTEMP=#1\relax
      \ifnum\dg@XTEMP<\@ne \dg@XTEMP=\@ne\fi
      \advance\dg@X\dg@XTEMP
      \ifnum\dg@HORIZ<\dg@X \dg@HORIZ=\dg@X \fi
      \setbox\dg@NODEBOX=\hbox{$\dgeverynode{#2}$}%
      \dg@XTEMP=\wd\dg@NODEBOX \advance\dg@XTEMP\dgHORIZPAD
      \dg@YTEMP=\ht\dg@NODEBOX \advance\dg@YTEMP\dp\dg@NODEBOX
      \advance\dg@YTEMP\dgVERTPAD
      \toks\z@=\expandafter{\dg@SLIST}%
      \edef\dg@SLIST{\the\toks\z@
         ,\noexpand\dg@XNODE=\number\dg@X\noexpand\relax
         \noexpand\dg@YNODE=\number\dg@Y\noexpand\relax
         \noexpand\dg@XTEMP=\number\dg@XTEMP\noexpand\relax
         \noexpand\dg@YTEMP=\number\dg@YTEMP\noexpand\relax}%
      \toks\z@=\expandafter{\dg@NLIST}%
      \toks\tw@={\dg@node{#2}}%
      \edef\dg@NLIST{\the\toks\z@
         \noexpand\dg@X=\number\dg@X\noexpand\relax
         \noexpand\dg@Y=\number\dg@Y\noexpand\relax
         \the\toks\tw@}%
   \else\ifx\dg@PASS\dg@geompass
      \ifnum\dg@X=\z@
         \dg@getnodesize
            {\dg@SLIST}{\dg@X}{\dg@Y}{\dg@WEND}{\dg@HEND}%
         \divide\dg@WEND\tw@
         \ifnum\dg@XLPAD<\dg@WEND \dg@XLPAD=\dg@WEND \fi\fi
      \ifnum\dg@X=\dg@HORIZ
         \dg@getnodesize
            {\dg@SLIST}{\dg@X}{\dg@Y}{\dg@WEND}{\dg@HEND}%
         \divide\dg@WEND\tw@
         \ifnum\dg@XRPAD<\dg@WEND \dg@XRPAD=\dg@WEND \fi\fi
      \ifnum\dg@Y=\z@
         \dg@getnodesize
            {\dg@SLIST}{\dg@X}{\dg@Y}{\dg@WEND}{\dg@HEND}%
         \divide\dg@HEND\tw@
         \ifnum\dg@YTPAD<\dg@HEND \dg@YTPAD=\dg@HEND \fi\fi
      \ifnum\dg@Y=-\dg@VERT\relax
         \dg@getnodesize
            {\dg@SLIST}{\dg@X}{\dg@Y}{\dg@WEND}{\dg@HEND}%
         \divide\dg@HEND\tw@
         \ifnum\dg@YBPAD<\dg@HEND \dg@YBPAD=\dg@HEND \fi\fi
   \else\ifx\dg@PASS\dg@drawpass
      \dg@XNODE=\dg@X \multiply\dg@XNODE\@m
      \multiply\dg@XNODE\dg@XGRID
      \dg@YNODE=\dg@Y \multiply\dg@YNODE\@m
      \multiply\dg@YNODE\dg@YGRID
      \setbox\dg@NODEBOX=\hbox{$\dgeverynode{#2}$}%
      \put(\dg@XNODE,\dg@YNODE){\dg@makebox{\box\dg@NODEBOX}}%
   \fi\fi\fi}%

\newoptcommand{\dg@cr}{\@ne}[1]{%
   \ifx\dg@PASS\dg@savepass
      \dg@YTEMP=#1\relax
      \ifnum\dg@YTEMP<\@ne \dg@YTEMP=\@ne \fi
      \advance\dg@Y -\dg@YTEMP\relax
      \dg@X=-\@ne\relax\fi}%
\newoptcommand{\dg@arrow}{\@ne}[2]{%
   \begingroup
   \dg@USERSIZE=#1\relax
   \ifnum\dg@USERSIZE<\@ne \dg@USERSIZE=\@ne \fi
   \dg@parse{#2}%
   \ifx\dg@PASS\dg@savepass
      \ifx\dg@label\dgl@b \let\dg@label=\dgl@t \fi
      \ifx\dg@label\dgl@r \let\dg@label=\dgl@l \fi
      \let\dg@process=\dg@save
   \else\ifx\dg@PASS\dg@geompass
      \let\dg@process=\dg@ignore
      \dg@geomcalc
   \else\ifx\dg@PASS\dg@drawpass
      \let\dg@process=\dg@draw
      \dg@drawcalc
   \fi\fi\fi
   \dg@label{\dg@process{#1}{#2}}}%

\newoptcommand{\arrow}{\@ne}[2]{%
   \dg@parse{#2}%
   \ifx\dg@label\dgl@b \let\dg@label=\dgl@t \fi
   \ifx\dg@label\dgl@r \let\dg@label=\dgl@l \fi
   \dg@label{\dg@textarrow{#1}{#2}}}%

\def\dg@textarrow#1#2#3#4{%
   \mathop{{\dgHORIZPAD=0pt\relax\dgVERTPAD=0pt\relax
      \begin{diagram}
         \node{}\arrow[#1]{#2}{#3}{#4}\node{}
      \end{diagram}}}}

\def\dg@parse#1{%
   \let\dg@label=\dgl@ \dgo@
   \let\dg@type=\@empty \let\dg@lbltype=\@empty
   \@for\dg@list:=#1\do{%
      \ifx\dg@type\@empty \let\dg@type=\dg@list
      \else\ifx\dg@lbltype\@empty \let\dg@lbltype=\dg@list
         \@ifundefined{dgo@\dg@list}{}{\@nameuse{dgo@\dg@list}}%
      \else
         \@ifundefined{dgo@\dg@list}{}{\@nameuse{dgo@\dg@list}}%
      \fi\fi}%
   \@ifundefined{dgt@\dg@type}{\dgt@e}{\@nameuse{dgt@\dg@type}}%
   \@ifundefined{dgl@\dg@lbltype}{}{%
      \dg@letname\dg@label{dgl@\dg@lbltype}}}

\def\dg@draw#1#2#3#4{%
   \put(\dg@X,\dg@Y){\dg@makebox{%
      \begin{picture}(0,0)%
         \thinlines
         \put(\dg@XOFFSET,\dg@YOFFSET){%
            \dg@VECTOR(\dg@DX,\dg@DY){\dg@SIZE}}%
         \put(\dg@XLBL,\dg@YLBL){\dg@makebox{%
            \begin{picture}(0,0)%
               \put(\dg@XLBLOFF,\dg@YLBLOFF){%
                  \dg@makebox[\dg@LBLONE]{$\dgeverylabel{#3}$}}%
               \put(-\dg@XLBLOFF,-\dg@YLBLOFF){%
                  \dg@makebox[\dg@LBLTWO]{$\dgeverylabel{#4}$}}%
            \end{picture}}}%
      \end{picture}}}%
   \endgroup}%

\def\dg@save#1#2#3#4{%
   \endgroup 
   \toks\z@=\expandafter{\dg@ALIST}%
   \toks\tw@={\dg@arrow[#1]{#2}{#3}{#4}}%
   \edef\dg@ALIST{\the\toks\z@%
      \noexpand\dg@X=\number\dg@X\noexpand\relax
      \noexpand\dg@Y=\number\dg@Y\noexpand\relax
      \the\toks\tw@}}%

\def\dg@ignore#1#2#3#4{\endgroup}

\def\dg@geomcalc{%
   \dg@XEND=\dg@SIZE \multiply\dg@XEND\dg@USERSIZE
   \ifnum\dg@DX=\z@
      \dg@YEND=\dg@XEND \dg@XEND=\z@
      \dg@changesign\dg@YEND\dg@DY
   \else
      \dg@changesign\dg@XEND\dg@DX \dg@YEND=\dg@XEND
      \multiply\dg@YEND\dg@DY \divide\dg@YEND\dg@DX
   \fi
   \advance\dg@XEND\dg@X \advance\dg@YEND\dg@Y
   \dg@getnodesize
      {\dg@SLIST}{\dg@XEND}{\dg@YEND}{\dg@WEND}{\dg@HEND}%
   \dg@XOFFSET=\dg@WEND \dg@YOFFSET=\dg@HEND
   \dg@getnodesize
      {\dg@SLIST}{\dg@X}{\dg@Y}{\dg@WEND}{\dg@HEND}%
   \advance\dg@XOFFSET\dg@WEND \divide\dg@XOFFSET\tw@
   \advance\dg@YOFFSET\dg@HEND \divide\dg@YOFFSET\tw@
   \dg@XTEMP=\dgARROWLENGTH \dg@YTEMP=\dgARROWLENGTH
   \ifnum\dg@DX>\z@
      \dg@ZTEMP=\dg@DX \multiply\dg@XTEMP\dg@DX
   \else \dg@ZTEMP=-\dg@DX \multiply\dg@XTEMP -\dg@DX \fi
   \ifnum\dg@DY>\z@
      \advance\dg@ZTEMP\dg@DY \multiply\dg@YTEMP\dg@DY
   \else \advance\dg@ZTEMP -\dg@DY \multiply\dg@YTEMP -\dg@DY\fi
   \ifnum\dg@ZTEMP=\z@\else
      \divide\dg@XTEMP\dg@ZTEMP \divide\dg@YTEMP\dg@ZTEMP
      \advance\dg@XOFFSET\dg@XTEMP \advance\dg@YOFFSET\dg@YTEMP
   \fi
   \divide\dg@XOFFSET\dg@SIZE \divide\dg@XOFFSET\dg@USERSIZE
   \divide\dg@YOFFSET\dg@SIZE \divide\dg@YOFFSET\dg@USERSIZE
   \ifnum\dg@DX=\z@ \dg@XOFFSET=\z@ \fi
   \ifnum\dg@DY=\z@ \dg@YOFFSET=\z@ \fi
   \ifnum\dg@XGRID<\dg@XOFFSET \global\dg@XGRID=\dg@XOFFSET\fi
   \ifnum\dg@YGRID<\dg@YOFFSET \global\dg@YGRID=\dg@YOFFSET\fi
   \relax}

\def\dg@drawcalc{%
   \dg@XEND=\dg@SIZE \multiply\dg@XEND\dg@USERSIZE
   \ifnum\dg@DX=\z@
      \dg@YEND=\dg@XEND \dg@XEND=\z@
      \dg@changesign\dg@YEND\dg@DY
   \else
      \dg@changesign\dg@XEND\dg@DX \dg@YEND=\dg@XEND
      \multiply\dg@YEND\dg@DY \divide\dg@YEND\dg@DX
   \fi
   \advance\dg@XEND\dg@X \advance\dg@YEND\dg@Y
   \dg@getnodesize
      {\dg@SLIST}{\dg@XEND}{\dg@YEND}{\dg@WEND}{\dg@HEND}%
   \divide\dg@WEND\unitlength \divide\dg@HEND\unitlength
   \multiply\dg@DX\dg@XGRID \multiply\dg@DY\dg@YGRID
   \dg@rmcommondiv\tw@\dg@DX\dg@DY
   \dg@rmcommondiv\tw@\dg@DX\dg@DY 
   \dg@rmcommondiv\thr@@\dg@DX\dg@DY
   \multiply\dg@SIZE\dg@USERSIZE \multiply\dg@SIZE\@m
   \ifnum\dg@DX=\z@
      \multiply\dg@SIZE\dg@YGRID
      \divide\dg@HEND\tw@ \advance\dg@SIZE -\dg@HEND
      \dg@getnodesize
         {\dg@SLIST}{\dg@X}{\dg@Y}{\dg@WEND}{\dg@YOFFSET}%
      \divide\dg@YOFFSET\unitlength \divide\dg@YOFFSET\tw@
      \advance\dg@SIZE -\dg@YOFFSET
      \dg@XOFFSET=\z@
      \def\dg@LBLONE{r}\def\dg@LBLTWO{l}%
      \dg@XLBL=\z@ \dg@YLBL=\dg@SIZE
      \multiply\dg@YLBL\dg@LBLPOS
      \divide\dg@YLBL\dgARROWPARTS\relax
      \advance\dg@YLBL\dg@YOFFSET
      \dg@changesign\dg@YLBL\dg@DY
      \dg@changesign\dg@YOFFSET\dg@DY
   \else
      \multiply\dg@SIZE\dg@XGRID
      \ifnum\dg@DY=\z@
         \divide\dg@WEND\tw@ \advance\dg@SIZE -\dg@WEND
         \dg@getnodesize
            {\dg@SLIST}{\dg@X}{\dg@Y}{\dg@XOFFSET}{\dg@HEND}%
         \divide\dg@XOFFSET\unitlength \divide\dg@XOFFSET\tw@
         \advance\dg@SIZE -\dg@XOFFSET
         \dg@YOFFSET=\z@
         \def\dg@LBLONE{b}\def\dg@LBLTWO{t}%
         \dg@YLBL=\z@ \dg@XLBL=\dg@SIZE
         \multiply\dg@XLBL\dg@LBLPOS
         \divide\dg@XLBL\dgARROWPARTS\relax
         \advance\dg@XLBL\dg@XOFFSET
         \dg@changesign\dg@XLBL\dg@DX
         \dg@changesign\dg@XOFFSET\dg@DX
      \else
         \divide\dg@WEND\tw@ \divide\dg@HEND\tw@
         \multiply\dg@HEND\dg@DX \divide\dg@HEND\dg@DY
         \ifnum\dg@HEND<\z@ \multiply\dg@HEND\m@ne \fi
         \ifnum\dg@WEND<\dg@HEND \advance\dg@SIZE -\dg@WEND
         \else \advance\dg@SIZE -\dg@HEND \fi
         \dg@getnodesize
            {\dg@SLIST}{\dg@X}{\dg@Y}{\dg@WEND}{\dg@HEND}%
         \divide\dg@WEND\unitlength \divide\dg@WEND\tw@
         \divide\dg@HEND\unitlength \divide\dg@HEND\tw@
         \multiply\dg@HEND\dg@DX \divide\dg@HEND\dg@DY
         \ifnum\dg@HEND<\z@ \multiply\dg@HEND\m@ne \fi
         \ifnum\dg@WEND<\dg@HEND \dg@XOFFSET=\dg@WEND
         \else \dg@XOFFSET=\dg@HEND \fi
         \advance\dg@SIZE -\dg@XOFFSET
         \dg@changesign\dg@XOFFSET\dg@DX
         \dg@YOFFSET=\dg@XOFFSET
         \multiply\dg@YOFFSET\dg@DY \divide\dg@YOFFSET\dg@DX
         \def\dg@LBLONE{br}\def\dg@LBLTWO{tl}%
         \ifnum\dg@DX<\z@ \ifnum\dg@DY>\z@
            \def\dg@LBLONE{bl}\def\dg@LBLTWO{tr}\fi\fi
         \ifnum\dg@DX>\z@ \ifnum\dg@DY<\z@
            \def\dg@LBLONE{bl}\def\dg@LBLTWO{tr}\fi\fi
         \dg@XLBL=\dg@SIZE
         \multiply\dg@XLBL\dg@LBLPOS
         \divide\dg@XLBL\dgARROWPARTS\relax
         \dg@changesign\dg@XLBL\dg@DX
         \dg@YLBL=\dg@XLBL
         \multiply\dg@YLBL\dg@DY \divide\dg@YLBL\dg@DX
         \advance\dg@XLBL\dg@XOFFSET
         \advance\dg@YLBL\dg@YOFFSET
      \fi
   \fi
   \dg@XLBLOFF=-\dg@DY \dg@changesign\dg@XLBLOFF\dg@DX
   \dg@YLBLOFF=\dg@DX \dg@changesign\dg@YLBLOFF\dg@DX
   \ifnum\dg@DX=\z@ \dg@XLBLOFF=\m@ne \fi
   \dg@XTEMP=\dg@DX \dg@changesign\dg@XTEMP\dg@DX
   \dg@YTEMP=\dg@DY \dg@changesign\dg@YTEMP\dg@DY
   \ifnum\dg@YTEMP>\dg@XTEMP \dg@XTEMP=\dg@YTEMP \fi
   \ifnum\dg@XTEMP=\z@ \dg@XTEMP=\@ne \fi
   \multiply\dg@XLBLOFF\dg@LBLOFF \divide\dg@XLBLOFF\dg@XTEMP
   \multiply\dg@YLBLOFF\dg@LBLOFF \divide\dg@YLBLOFF\dg@XTEMP
   \multiply\dg@X\@m \multiply\dg@X\dg@XGRID
   \multiply\dg@Y\@m \multiply\dg@Y\dg@YGRID
   \relax}%

\def\dg@rmcommondiv#1#2#3{%
   \dg@XTEMP=#2\relax
   \divide\dg@XTEMP #1\relax \multiply\dg@XTEMP #1\relax
   \dg@YTEMP=#3\relax
   \divide\dg@YTEMP #1\relax \multiply\dg@YTEMP #1\relax
   \ifnum\dg@XTEMP=#2\relax \ifnum\dg@YTEMP=#3\relax
      \divide#2#1\relax \divide#3#1\relax \fi\fi}%

\def\dg@changesign#1#2{%
   \ifnum #2<\z@ \multiply#1\m@ne
   \else\ifnum #2=\z@ #1=\z@ \fi\fi}%

\def\dg@getnodesize#1#2#3#4#5{%
   #4=\z@\relax #5=\z@\relax
   \expandafter\@for\expandafter\dg@trynode
   \expandafter:\expandafter=#1\do{%
      \dg@XNODE=\m@ne 
      \dg@trynode
      \ifnum #2=\dg@XNODE \ifnum #3=\dg@YNODE
         #4=\dg@XTEMP\relax #5=\dg@YTEMP\relax\fi\fi}}%

\newoptcommand{\dg@makebox}{}[2]{%
   \expandafter\makebox\expandafter(\expandafter
      0\expandafter,\expandafter0\expandafter)\expandafter
      [#1]{#2}}%

\def\dg@novector(#1,#2)#3{}%

\def\dg@letname#1#2{%
   \relax\expandafter
   \let\expandafter #1\csname #2\endcsname\relax}%

\def\dgl@#1{#1{}{}}%
\def\dgl@t#1#2{#1{#2}{}}%
\def\dgl@b#1#2{#1{}{#2}}%
\def\dgl@tb#1#2#3{#1{#2}{#3}}%
\def\dgl@l#1#2{#1{#2}{}}%
\def\dgl@r#1#2{#1{}{#2}}%
\def\dgl@lr#1#2#3{#1{#2}{#3}}%
\makeatother

\theoremstyle{plain} 

\newtheorem{thm}{Theorem}
\newtheorem{theo}{Theorem}[section]
\newtheorem{lem}[theo]{Lemma}
\newtheorem{cor}[theo]{Corollary}
\newtheorem{prop}[theo]{Proposition}

{\theorembodyfont{\rmfamily} \newtheorem{rem}[theo]{Remark}}
{\theorembodyfont{\rmfamily} }

\newcommand{\op}[1]{\operatorname{#1}}

\newcommand{\bbp}{{\mathbb P}}
\newcommand{\calo}{{\mathcal O}}

\newcounter{dummy}

\begin{document}
\title{\bf Weak Hironaka Theorem}
\author{Fedor A. Bogomolov\thanks{Partially supported by NSF grant 
DMS-9500774.} \and Tony G. Pantev\thanks{Partially supported by NSF grant 
DMS-9500712.}}
\date{ }
\maketitle

\tableofcontents

\section{Introduction} \label{s1}

The existence of a  smooth projective model for any proper algebraic
variety  over an algebraically closed field of characteristic zero is  one 
of the most important results in algebraic geometry.

H.Hironaka \cite{h} proved a very strong version of the above statement. 
For any such variety he established the existence of a sequence 
of blow ups which resolve the singularities of the variety. Moreover at every 
step  the blow ups occur along smooth subvarieties of the singular locus. 
In particular the process does not modify the nonsingular  part of the initial 
variety. This result proved to be even more useful than the existence
of a smooth model, but its proof was rather complicated.
Later on, E. Bierstone and P. Milman \cite{m}
and M. Spivakovsky \cite{s} found different versions of the resolution process
that are considerably simpler and canonical but somehow in principle their 
proofs follow the same track.

In this note we want to give a proof 
of the existence of a smooth projective model
for any variety in characteristic zero using an idea similar to A.J. de Jong's 
approach in \cite{dj}. 
We present a  straightforward way to obtain a nonsingular model which 
however does not permit any control over the intermediate steps.   
 
The main result of this note is the following theorem

\begin{thm} \label{thm1} Let $X$ be a normal projective variety 
over an algebraically closed field $k$, $\op{char} k = 0$ and let 
$D \subset X$ be a 
proper subvariety of $X$. Then there exist a smooth projective variety
$M$, a strict normal crossings  divisor $R \subset M$  and a birational 
morphism $f : M \to X$ with $f^{-1} D = R$
\end{thm}

Our proof uses induction on the dimension of $X$. It is organized as 
follows. In section~2  we construct a finite 
morphism $(\widehat{X},\widehat{D}) \to (P,S_{0})$ from a suitable blow up 
$(\widehat{X},\widehat{D})$ of the pair $(X,D)$ to a pair $(P,S_{0})$ where 
$P$ is a compactification of the total space of a line bundle over 
$\bbp^{n-1}$, $S_{0}$
is the image of a section of this line bundle and $\widehat{X} \to P$
is branched over a divisor $B = \sum S_{i}$ which is a sum of distinct 
sections of $P \to \bbp^{n-1}$. 
After this is achieved we proceed in section~3 by applying the induction
hypothesis to the pair $(\bbp^{n-1},Z)$. Here $Z \subset \bbp^{n-1}$ is the
discriminant locus of the set of sections $\{ S_{i} \}$, i.e. $Z$ is the 
unique closed reduced subscheme in $\bbp^{n-1}$ whose preimage in $P$ contains
all possible intersections of the $S_{i}$'s.  
By induction we can find a birational morphism $\widetilde{\bbp}^{n-1}
\to \bbp^{n-1}$ which transforms $Z$
into a strict normal crossings divisor $\widetilde{Z}$. Put $\widetilde{P}$ for
the variety obtained from $P$ via a  base change with $\widetilde{\bbp}^{
n-1}\to \bbp^{n-1}$. Then $(\widetilde{P} \to \widetilde{\bbp}^{n-1}, 
\{S_{i}\})$ is a family of pointed $\bbp^{1}$'s and an application of 
F. Knudsen's stabilization theorem yields a smooth variety
$Q$ equipped with a contraction map  onto $\widetilde{P}$ and such that the 
complete preimage of $B$ in $Q$ is a normal crossings divisor. Now the fiber
product of $Q$ and $\widehat{X}$ has only abelian quotient singularities 
which can be resolved by a toroidal blow-up. The resulting pair (M,R) will be 
the resolution of $(X,D)$. 

\begin{rem} \label{rem11} Recently A.J. de Jong \cite{dj} invented a 
remarkable new approach which allowed him to 
resolve in any characteristic at least some variety which dominates 
the initial variety by a finite map.
The idea to fiber $X$ by curves 
and use induction on the dimension came from de Jong's insight in 
\cite{dj}. The important reduction in section~2 imitates the first part of
the proof of the famous theorem of Belyi \cite{b}.
\end{rem}

\begin{rem} \label{rem12} The fact that $\op{char} k = 0$ is used only 
at last step of the argument in section~3 since in the case of a positive 
characteristic it is hard to control the unramified coverings of the torus.
\end{rem}
\begin{rem} \label{rem13} A different proof of the same theorem was found 
independently by {\em D. Abramovich and A.J. de Jong} \cite{A-dJ}. It uses 
similar ingredients but is quite different in detail. In particular their 
proof has the advantage of being better adapted to handling equivariant 
resolutions.
\end{rem}

\bigskip

\noindent
{\bf Acknowledgments} We would like to thank A.J. de Jong for his inspiring 
Santa Cruz lectures and D. Abramovich for suggesting that we use semistable
reduction for pointed curves of genus zero which allowed us to substantially
simplify our original application of the induction. We would also like to 
thank the referee for the careful reading of the manuscript and for the many
valuable comments.

\section{The map to a $\bbp^1$ bundle}

Let $X$ be an irreducible variety of dimension $n$ over $k$ and let $D \subset
X$ be a closed subvariety.  
The aim of this section is to prove the following theorem.

\begin{theo} \label{theo21}
There exists a finite 
morphism $(\widehat{X},\widehat{D}) \to (P,S_{0})$ from a suitable blow up 
$(\widehat{X},\widehat{D})$ of the pair $(X,D)$ to a pair $(P,S_{0})$ of 
varieties satisfying:

\begin{list}{{\em (\roman{dummy})}}{\usecounter{dummy}}
\item $P = \bbp(\calo\oplus L)$ for some line bundle $L \to \bbp^{n-1}$.
\item $S_{0} \subset P$ is the image of a section of $L$.
\item The finite morphism $\widehat{X} \to P$ is branched over a divisor
$B$ of the form 
\[
B = \sum_{i=1}^{k} S_{i} + S_{k+1},
\]
where the $S_{i}$'s, $i= 1,\ldots, k$ are images of distinct sections of 
$L$ and $S_{k+1} := \bbp_{\infty}$ is the infinity section of 
$P \to \bbp_{n-1}$.
\end{list}
\end{theo}
{\bf Proof.} Passing to a  blow up of $(X,D)$ if necessary we may 
assume without loss of generality that $X$ is normal and that 
$D = \sum D_i$ with $D_{i}$ reduced and irreducible divisors.
We start with the following relative version of the Noether normalization 
lemma.

\begin{lem}There exists a finite morphism $f : X \to \bbp^n$ which maps
$\widetilde{D}$ onto a hyperplane $\bbp^{n-1}$ of $\bbp^n$. \label{lem21}
\end{lem}
{\bf Proof.} In order to find such a map it is enough to consider the
embedding of $X$ into projective space via a very ample line bundle
which has a section vanishing on every $D_i$. This can be easily
achieved by taking a sufficiently big line bundle.
Let $X\subset \bbp^{m}$ be such an embedding and let $H \subset \bbp^{m}$
be the hyperplane containing all $D_i$. A straightforward dimension count
shows that the generic 
subspace $\bbp^{N-n-1}\subset H$ does not intersect $X\cap H$ and therefore
we can take for $f$ the projection of $\bbp^m$ centered at such a subspace.
\hfill $\Box$

\begin{lem} The map $f$ can be chosen so that
there exists a point $o \in \bbp^n$ satisfying the following properties:
\begin{list}{{\em (\alph{dummy})}}{\usecounter{dummy}}
\item the preimage of $o$ consists
of a finite number of smooth points $x_i$ with $df_{x_i}$ 
 being an isomorphism for every $x_{i}$.
\item the preimage of every line in $\bbp^n$ containing $o$ is a connected
and generically reduced curve. 
\item the preimage of the generic line trough $o$ is irreducible.
\end{list} \label{lem22}
\end{lem}
{\bf Proof.} The subset of points satisfying (a) is open
in $\bbp^n$ due to the generic smoothness lemma. In order
to satisfy (b) it suffices to chose the projective embedding from 
lemma~\ref{lem21} in such a way that all non-reduced curves in $X$
obtained as hyperplane sections constitute a subvariety  of big codimension
(greater than $n-1$) in the corresponding Grassmanian. In other words we 
want to choose the embedding $X \subset \bbp^{m}$ so that there exists a
projective subspace $V \subset \bbp^m$ of dimension $m-n$ with the property 

\medskip

\noindent
\begin{description}
\item[(*)] For every projective subspace $V \subset M \subset \bbp^m$ of 
dimension $m-n+1$ the intersection $M\cap X$ is a generically reduced 
connected curve.
\end{description}

\medskip

\noindent
Due to the connectedness theorem \cite[III 
Corollary 7.9]{ha}, \cite{fh} the preimage in $X$ of every line trough $o$
will be connected. Let $Y$ be the union of the singular locus of $X$ the 
divisor $D$ and the ramification divisor of $f$. Furthermore 
Bertini's theorem implies that the part of a general  
hyperplane curve contained in $X\setminus Y$ will be reduced and irreducible 
and therefore it suffices to show that for a suitable embedding of $X$ there 
exists a $m-n$ dimensional projective subspace $V \subset  \bbp^{m}$
satisfying

\medskip

\noindent
\begin{description}
\item[(**)] For every projective subspace $V \subset M \subset \bbp^m$
of dimension $m-n+1$ the intersection $M\cap Y$ is of dimension zero.
\end{description}

\medskip
\noindent
Consider the dimension jump locus 
\[
\Gamma = \left\{ M \left| \dim M = m-n+1, \; \dim(M\cap Y) \geq 
1 \right. \right\} \subset Gr(m-n+2,m+1),
\]
and let $\widetilde{\Gamma} = \{ (M,y) | y \in M \} \subset \Gamma \times Y$ 
be the corresponding incidence variety. The projection $\pi_{Y} : 
\widetilde{\Gamma} \to Y$ surjects on $Y$ and its fiber over a general 
point $y \in Y$ consists of all $M \subset \bbp^m$ of dimension $m-n+1$
such that $y \in M$ and $\dim(M\cap Y) \geq 1$. In particular $\dim 
\widetilde{\Gamma} \leq 2n - 3$ and $\dim \Gamma \leq 2n -4$. Furthermore
the flag variety $F$ parametrizing flags $V \subset M \subset \bbp^m$ as
above has the usual double fibration structure
\[
{\divide\dgARROWLENGTH by 5
\begin{diagram}
\node[2]{F} \arrow{sw,l}{p} \arrow{se,l}{q} \\
\node{Gr(m-n+1,m+1)} \node[2]{Gr(m-n+2,m+1)}
\end{diagram}}
\] 
where the fibers of $p$ are projective spaces of dimension $n-1$ and the fibers
of $q$ are projective spaces of dimension $m-n$. The locus in $Gr(m-n+1,m+1)$
consisting of all $V$'s that do not satisfy (**) is contained in
$q(p^{-1}(\Gamma))$ and hence has codimension that is $\geq (m-n+1)n - (m-n + 
2n -4) = m(n+1) - n^2 + 4$. Thus composing our original embedding $X \subset 
\bbp^m$ with a Veronese embedding of sufficiently high degree we can make this
codimension strictly positive which proves the lemma. \hfill $\Box$

\bigskip

Denote by $P_{1}$ the blow up of $\bbp^{n}$ at $o$ and by $\widehat{X}$ the
blow up of $X$ at the points $x_i$. Since the blow up at $o \in \bbp^n$
resolves the projection centered at $o$ we get a realization of $P_{1}$ as a 
$\bbp^1$ bundle over $\bbp^{n-1}$. There is a canonical isomorphism
$P_{1} \cong \bbp(\calo_{\bbp^{n-1}}\oplus \calo_{\bbp^{n-1}}(1))$ which 
gives a natural tautological line bundle $\calo_{P_{1}}(1) \to P_{1}$ and the
preimage $\bbp_{\infty}$  of $o$ in $P_1$ is a section of $\calo_{P_{1}}(1)$.
Alternatively, $\bbp_{\infty}$ is the divisor in $\bbp(\calo_{\bbp^{n-1}}
\oplus \calo_{\bbp^{n-1}}(1))$ corresponding to the line subbundle 
$\calo_{\bbp^{n-1}} \subset \calo_{\bbp^{n-1}}\oplus \calo_{\bbp^{n-1}}(1)$
and $P_1$ can be thought of as a compactification of the total space of 
$\calo_{\bbp^{n-1}}(1)$ with  the divisor $\bbp_{\infty}$ at infinity.

By construction  $\widehat{X}$ is fibered by connected and generically reduced
curves $C_t, t\in \bbp^{n-1}$ and 
maps fiberwise via $f$ to the $\bbp^1$-fibration $P_1$.
Due to the smoothness of $P_{1}$ the Zariski-Nagata purity theorem 
\cite[X 3.1]{sga} implies that the branch locus $B_{1}$ of the
finite map  $f : \widehat{X} \to P_1$ is a divisor. Moreover, since $C_{t}$ is 
generically reduced for every $t$ we conclude that $B_{1}$ is a horizontal
 divisor in $\pi : P_{1} \to \bbp^{n-1}$,
i.e. every component of $B_{1}$ is mapped finitelyonto $\bbp^{n-1}$ by $\pi$.
Denote by $N$ the degree of $B_{1}$ on the fibers of $\pi$.
Now we can apply Belyi's construction to simplify the branch divisor.

We will need the following technical result:

\begin{lem} Let $L$ be a line bundle on a variety $M$ and $B$ be a horizontal
divisor in the total space $\op{tot}(L)$ of $L$, which maps properly on $M$ 
under a natural projection $\pi : \op{tot}(L) \to M$.
Denote by $d$ the degree of the finite map $\pi : B \to M$.
There exists a canonical fiberwise morphism $p_B : \op{tot}(L) \to 
\op{tot}(L^{\otimes d})$ which is polynomial of degree $d$ on any fiber $L_x, 
x\in M$ and such that $B$ is the preimage of the zero section of $L^{\otimes 
d}$ under $p_B$. The map $p_B$ is defined uniquely modulo multiplication by 
invertible functions on $M$.
\end{lem}
{\bf Proof.} Consider the compactification $P_{L} := \bbp(\calo_M \oplus L) =
\op{Proj}(\op{Symm}^{\bullet}(\calo_{M}\oplus L^{-1}))$.
Let $M_{\infty}$ be the infinity section of $P_L$ and $y \in H^{0}(M,
\mathcal{O}_{M}) \subset H^{0}(P_{L}, \calo_{P_{L}}(1))$ be a non-zero 
element. Then $\op{div}(y) = M_{\infty}$.
Denote by $\lambda \in H^{0}(\op{tot}(L), \pi^*L)$ the tautological section
and let $x \in H^{0}(P_{L}, \calo_{P_{L}}(1)\otimes \pi^{*}L)$ be the 
section whose divisor is the zero section of $L$ normalized so that 
$\lambda = x/y$. The fact that the intersection of $B$ with the generic 
fiber of $\pi : P_L \to M$ is $d$ implies that $\calo_{P_L}(B) = 
\calo_{P_{L}}(d)\otimes\pi^*A$ for some line bundle $A$ on $M$. 
Moreover since $B \subset \op{tot}(L)$
it follows that $\calo(B)_{|M_{\infty}} = \calo$ and hence $\pi^{*}
A_{|M_{\infty}} \cong \calo_{P_L}(-d)_{|M_{\infty}}$.  On the
other hand $P_L =  \bbp(\calo_M \oplus L)$ and by using that 
$\pi_{|M_{\infty}}$ is an isomorphism we get that $A = L^{\otimes d}$.

Let $\varphi
\in H^{0}(P_{L},\calo_{P_{L}}(d)\otimes\pi^*L^{\otimes d})$ be a section 
with divisor $B$. Applying the projection formula to the pushforward by $\pi$ 
we get isomorphisms
\[
\begin{aligned}
H^{0}(P_{L}, & \calo_{P_{L}}(d)\otimes\pi^*L^{\otimes d}) = 
H^{0}(M,(\pi_{*}\calo_{P_{L}}(d))\otimes L^{\otimes d}) = \\  
& = H^{0}(M,S^d(\calo_M\oplus L^{-1})\otimes L^{\otimes d}) = \\
& = H^{0}(M,\calo_M)\oplus H^{0}(M,L) \oplus \ldots \oplus 
H^{0}(M,L^{\otimes d}).
\end{aligned}
\]
The image of $\varphi$ under this isomorphism can be decomposed as 
$\varphi = (\varphi_{0}, \ldots, \linebreak \varphi_{d})$ with $\varphi_{i} 
\in H^{0}(M,L^{\otimes i})$ and hence
\begin{equation} \label{eq21}
\varphi = (\pi^*\varphi_{0})x^d +(\pi^*\varphi_{1})x^{d-1}y + \ldots 
+ (\pi^*\varphi_{d})y^d, 
\end{equation}
and
\begin{equation} \label{eq22}
\varphi_{|\op{tot}(L)} = (\pi^*\varphi_{0})\lambda^d +(\pi^*\varphi_{1})
\lambda^{d-1} + \ldots 
+ (\pi^*\varphi_{d}).
\end{equation}
Now the map 
\[
\begin{array}{llcl}
p_{B} : & \op{tot}(L) & \longrightarrow & \op{tot}(L^{\otimes d}) \\
 & \eta & \longrightarrow & (\pi^*\varphi_{0})\eta^d +(\pi^*\varphi_{1})
\eta^{d-1} + \ldots 
+ (\pi^*\varphi_{d})
\end{array}
\]
has the desired properties. \hfill $\Box$

\begin{cor} The above map extends to a proper map of the compactifications
$p_B :  P_L \to P_{L^{\otimes d}}$. The preimage of infinite section 
$M_{\infty, L^{\otimes d}}$ is a $d$-multiple of $M_{\infty, L}$ \label{cor21}
\end{cor}
{\bf Proof.} The extension of $p_B$ is given by the right hand side of 
(\ref{eq21}). \hfill $\Box$

\begin{cor} The map $p_B :  \op{tot}(L) \to \op{tot}(L^{\otimes d})$ is 
branched over a
divisor $B' \subset \op{tot}(L^{\otimes d})$ which is horizontal and of degree 
$d-1$ along the fibers. \label{cor22}
\end{cor}
{\bf Proof.} From the definition of $p_{B}$ it is clear that the divisor 
$B'$ is just the image under $p_{B}$ of the divisor of the section
$p_{B}'(\lambda) := d(\pi^*\varphi_{0})\lambda^{d-1} + \linebreak 
(d-1)(\pi^*\varphi_{1})
\lambda^{d-2} + \ldots + (\pi^*\varphi_{d-1}) \in H^{0}(\op{tot}(L),
\pi^*L^{\otimes (d-1)})$. \hfill $\Box$

\begin{rem} Corollary~\ref{cor22} is a complete analogue of the crucial 
observation of Belyi in \cite{b}.
\end{rem}

Now we can finish the proof of theorem~\ref{theo21}. Indeed let $P_{2} = 
\bbp(\calo_{\bbp^{n-1}}\oplus \calo_{\bbp^{n-1}}(N))$ and let $p_{B_{1}} : 
P_{1} \to P_{2}$ be as in corollary~\ref{cor21}. Then according to 
corollary~\ref{cor22} the branch divisor of $p_{B_{1}}\circ f$ is the union 
of the zero section of $\calo_{\bbp^{n-1}}(N)$ and a horizontal divisor $B_2$
of degree $N-1$ along the fibers. Now we can compose with $p_{B_{2}} + s$
where $s$ is a generic section of $\calo_{\bbp^{n-1}}(N(N-1))$ and continue
until we get a map to a $\bbp^1$ bundle over $\bbp^{n-1}$ branched over a 
union of distinct sections. \hfill $\Box$

\section{The induction step} In this section we use the induction hypothesis 
to modify birationally the pair $(P,B\cup S_{0})$ obtained in 
Theorem~\ref{theo21} so that $B\cup S_{0}$ becomes a divisor with strict 
normal crossings. First we need the following proposition.

\begin{prop} \label{prop31} There exists a diagram
\[
{\divide\dgARROWLENGTH by 5
\begin{diagram}
\node{Q} \arrow{e,t}{\varepsilon} \arrow{se,b}{\pi} \node{\widetilde{P}} 
\arrow{e} \arrow{s,l}{\tilde{p}} \node{P} \arrow{s,r}{p} \\
\node[2]{\widetilde{\bbp}^{n-1}} \arrow{e} \node{\bbp^{n-1}} \\
\node[2]{\widetilde{Z}} \arrow{n,J} \arrow{e} \node{Z} \arrow{n,J} 
\end{diagram}}
\]
where
\begin{list}{{\em (\roman{dummy})}}{\usecounter{dummy}}
\item $\widetilde{\bbp}^{n-1} \to \bbp^{n-1}$ is a birational morphism. 
$\widetilde{P} = \widetilde{\bbp}^{n-1}\times_{\bbp^{n-1}} P$ and 
$\widetilde{Z}=\widetilde{\bbp}^{n-1} \times_{\bbp^{n-1}} Z$;
\item $\pi$ is a flat and proper morphism whose geometric fibers are
reduced and connected curves with at most ordinary double points. 
Furthermore $\varepsilon : Q\setminus \pi^{-1}(\widetilde{Z}) \to \widetilde{P}
\setminus\tilde{p}^{-1}(\widetilde{Z})$ is an isomorphism;
\item $\varepsilon$ is a birational morphism which restricted 
on a geometric fiber $Q_{x}$ of $\pi$ falls under one of the following 
possibilities:
\begin{list}{{\em \alph{dummy})}}{\usecounter{dummy}}
\item $\varepsilon_{x} : Q_{x} \to \widetilde{P}_{x}$ is an isomorphism;
\item There is a (possibly disconnected) reduced rational curve $E \subset 
Q_{x}$ having at most ordinary double points so that $\varepsilon_{x}(E) = 
\Sigma$ is a finite set of closed points in $\widetilde{P}_{x}$ and
\[
\varepsilon_{x} : Q_{x}\setminus E \to \widetilde{P}_{x}\setminus \Sigma
\]
is an isomorphism;
\end{list} 

\item $\widetilde{Z}$ and the preimage of $B\cup S_{0}$ in $Q$ 
are strict normal crossings divisors;
\item $Q$ is non-singular.
\end{list} 
\end{prop}
{\bf Proof.} Denote by $Z\subset \bbp^{n-1}$ the discriminant locus of 
the set of sections $S_{0}, S_{1}, \ldots, S_{k+1}$. That is $Z$ is the 
smallest closed reduced subscheme in $\bbp^{n-1}$ whose preimage in $P$ 
contains all possible intersections of the $S_{i}$'s. By induction there 
exists a birational morphism $\mu: \widetilde{\bbp}^{n-1} \to \bbp^{n-1}$  
such that the reduced preimage $\widetilde{Z}$ of $Z$ is a divisor with strict 
normal crossings. Denote by $\widetilde{P}$ the fiber product of $P$ with 
$\widetilde{\bbp}^{n-1}$. Clearly $\widetilde{P} = \bbp(\mu^{*}L\oplus 
\mathcal{O})$. Denote by $\tilde{s}_{i}$, $i = 0,\ldots, k+1$ the sections of 
the projective bundle $\tilde{p} : \widetilde{P} \to \widetilde{\bbp}^{n-1}$
whose images are the divisors $\widetilde{S}_{i} := \mu^{*}S_{i}$, $i = 0, 
\ldots ,k+1$.

The image $\widetilde{S}_{k+1}$ of the infinity section of  $\widetilde{p}$
does not intersect any of the divisors $\widetilde{S}_{i}$, $i = 0,\ldots, k$
and thus will not cause any trouble. To deal with the bad intersections of 
the remaining $k+1$ sections we appeal to F. Knudsen's stabilization theorem
\cite[Theorem~2.4]{k}. According to this theorem
given a flat family $C \to S$
of connected reduced curves with at most ordinary double points and $a$ 
distinct sections $\{s_{i} : S \to C\}_{i=1}^{a}$ plus an arbitrary extra 
section $\Delta : S \to C$ there exist a canonical blow-up $q: C' \to C$ and 
unique liftings of the sections $s_{1}, \ldots, s_{a}$ and $\Delta$ to sections
$s_{1}', \ldots, s_{a}',s_{a+1}'$ of $C' \to S$ so that for the induced 
morphism $q_{s} : C'_{s} \to C_{s}$ on any geometric fiber $C'_{s}$ of 
$C' \to S$ we have one of the following two cases 
\begin{list}{\alph{dummy})}{\usecounter{dummy}}
\item $q_{s} : C'_{s} \to C_{s}$ is an isomorphism;
\item There is a rational component $E \subset C'_{s}$ such that $s_{a+1}(s)
\in E$, $q_{s}(E) = x$ is a closed point of $C_{s}$ and
\[
q_{s} : C'_{s}\setminus E \longrightarrow C_{s}\setminus \{ x\} 
\]
is an isomorphism.
\end{list} 

\medskip

\noindent
In particular, $C' \to S$ is again a flat family of reduced curves with
at most ordinary double points. Explicitly $q: C' \to C$ is given as
follows. Denote by $\mathcal{J}$ the $\mathcal{O}_{C}$ ideal defining
$\Delta$ and by $\mathcal{C}$ the cokernel of the diagonal embedding
\[
\mathcal{O}_{C} \longrightarrow \mathcal{J}^{\vee}\oplus \mathcal{O}_{C}(s_{1}
+\ldots + s_{a}).
\] 
Then 
\[
C' := \op{Proj}(\oplus_{i \geq 0} \mathcal{C}^{i}),
\]
and $q : C' \to C$ is the natural structure morphism.

The total space of $C'$ can be made smooth by a sequence of blow-ups with 
smooth centers as long as $C$ 
and $S$ are smooth and the morphism $C \to S$ is smooth in the complement 
of a divisor with strict normal crossings. This can be seen directly by 
analyzing the local picture around the separated sections but instead of 
doing that we prefer to invoke a much  stronger general result of A.J. de 
Jong.  The \cite[Proposition~5.6]{dj} 
guarantees that the total space of any split flat family of semistable 
curves over a smooth base which is smooth over the complement of a strict 
normal crossings divisor can be blown-up so that the resulting family is of the
same type and  has a smooth total space. Thus the only thing that requires 
checking is that the family $C'$ is split but this is clear since we can 
label the components of every fiber of $C'$ by their distance from the 
proper transform of the fiber in $C$.

\begin{rem}\label{rem31} The stabilization theorem is stated in \cite{k} 
for families of stable curves since Knudsen is interested in obtaining a 
morphism between the universal curve over the stack of $n$-pointed curves and
the stack of $n+1$ pointed curves. However the proof he gives in \cite{k}
is carefully designed to work for general flat families of reduced nodal 
curves. In particular, the stability assumption is never used in his argument. 
\end{rem} 

\medskip

\noindent
Applying the stabilization theorem to each of the sections $\widetilde{S}_{1},
\ldots \widetilde{S}_{k}$ (or rather to the corresponding proper transforms) 
in turn, we obtain a variety $Q$ together with a  birational morphism 
$\varepsilon : Q \to \widetilde{P}$ and $k+2$ distinct sections $t_{i} : 
\widetilde{\bbp}^{n-1} \to Q$, $i = 0, \ldots,k+1$  lifting the sections 
$\tilde{s}_{i}$. Since at every step we are blowing-up
only ideal sheaves whose support is contained in the preimage of 
$\widetilde{Z}$ the construction of $Q$ guarantees the validity of items
(ii) and (iii) of the statement of the proposition. Furthermore due to  part 
b) of the stabilization theorem  the images of the sections $t_{i}$ intersect 
every fiber of $\pi$ at a smooth point. Thus the preimage of 
$B\cup S_{0}$ in $Q$ is the same as the union of the images of the $t_{i}$'s 
and $\pi^{-1}(\widetilde{Z})$. Since $\pi$ is flat and has singular fibers that
are trees of rational curves we conclude that $\pi^{-1}(\widetilde{Z})$ is 
also a divisor with strict normal crossings which finishes the proof of the
proposition. \hfill $\Box$

\bigskip

To finish the proof of Theorem~A consider the 
normalization $Y$ of the fiber product $X\times_{P} Q$. Denote by $Y^{\circ}$ 
the preimage of $Q\setminus \{ B 
\cup S_{0} \}$ in $Y$. The normal variety $Y$ is a finite cover of the
smooth variety $Q$ branched along a divisor with strict normal crossings.
Since we are in characteristic zero this implies that $Y$ has only abelian 
quotient singularities and is toroidal without self-intersections, i.e. 
locally in the \'{e}tale topology (or formally) the embedding $Y^{\circ} 
\subset Y$ is isomorphic to the embedding of an affine algebraic torus 
$U$ in its toric variety $\overline{U}$ so that every component of the 
complement $\overline{U}\setminus U$ is normal. Finally we invoke the 
fundamental theorem \cite[Ch. II, \S 2, Theorem~11$\ast$]{kkms} about toric 
resolutions according to which there exists a canonical sheaf of ideals on $I 
\subset \mathcal{O}_{Y}$ such that the blow-up $Bl_{I}(Y)$ is non-singular.

\bigskip

\noindent
{\sc Fedor Bogomolov \\
Courant Institute for Mathematical Sciences \\
New York University \\
New York, NY 10012-1110 \\
e-mail: bogomolo@@cims.nyu.edu}
 
\smallskip

{\sc and }

\smallskip

\noindent
{\sc Steklov Institute of Mathematics \\
Russian Academy of Sciences\\
42 Vavilova street \\
Moscow 117333, Russia}

\bigskip 

\bigskip

\bigskip

\noindent
{\sc  Tony Pantev \\
Department of Mathematics \\
Massachusetts Institute of Technology \\
Cambridge, MA 02139-4307 \\
e-mail: pantev@@math.mit.edu}

\end{document}